\documentstyle[12pt]{article}

\newfont{\kh}{msbm10}
\newcommand{\R}{\mbox{\kh R}}
\newcommand{\C}{\mbox{\kh C}}
\newcommand{\T}{\mbox{\kh T}}
\newtheorem{t.}{Theorem}[section]
\newtheorem{d.}[t.]{Definition}
\newtheorem{l.}[t.]{Lemma}
\newtheorem{p.}[t.]{Proposition}
\newtheorem{c.}[t.]{Corollary}
\newtheorem{e.}[t.]{Example}
\newtheorem{r.}[t.]{Remark}
\begin{document}
\title {Approximately Vanishing of Topological Cohomology Groups }
\author{M. S. Moslehian}
\date{}
\maketitle

\begin{abstract}
In this paper, we establish the Pexiderized stability of coboundaries and cocycles and use them to investigate the Hyers--Ulam stability of some functional equations. We prove that for each Banach algebra $A$, Banach $A$-bimodule $X$ and positive integer $n, H^n(A,X)=0$ if and only if the $n$-th cohomology group approximately vanishes.
\end{abstract}

\section{Introduction.}

Topological cohomology arose from the problems concerning extensions by H. Kamowitz who introduced the Banach version of Hochschild cohomology groups in 1962 $\cite{KAM}$, derivations by R. V. Kadison and J. R. Ringrose $\cite{K-R1}, \cite{K-R2}$ and amenability by B.E. Johnson $\cite{JOH}$ and has been extensively developed by A. Ya. Helemskii and his school $\cite{HEL1}$. In addition, this area includes a lot of problems concerning automorphism groups of operator algebras, fixed point theorems, stability, perturbations, invariant means $\cite{HEL1}$ and their applications to quantum physics $\cite{SAK}$.

Consider the functional equation ${\cal E}_1(f)={\cal E}_2(f)~~~({\cal E})$ in a certain framework. We say a function $f_0$ is an approximate solution of $({\cal E})$ if ${\cal E}_1(f_0)$ and ${\cal E}_2(f_0)$ are close in some sense. The stability problem is whether or not there is a true solution of $({\cal E})$ near $f_0$.

The stability of functional equations started with the following question concerning stability of group homomorphisms proposed by S. M. Ulam during a talk before a Mathematical Colloquium at the University of Wisconsin, Madison, in 1940: 

{\small Let $G_1$ be a group and let $(G_2, d)$ be a metric group. Given $\epsilon>0$, does there exist a $\delta>0$ such that if a mapping $f : G_1 \to G_2$ satisfies the inequality $d(f(xy), f(x)f(y))<\delta$ for all $x, y \in G_1$ then a homomorphism $T:G_1\to G_2$ exists such that $d(f(x),T(x))<\epsilon$ for all $x\in G_1$.}

In 1941, D. H. Hyers $\cite{HYE}$ provide the first (partial) answer to Ulam's problem as follows:

{If $E_1, E_2$ are Banach spaces and $f:E_1 \to E_2$ is a mapping for which there is $\epsilon>0$ such that $\|f(x+y)-f(x)-f(y)\|<\epsilon$ for all $x, y\in E_1$, then there is a unique additive mapping $T:E_1\to E_2$ such that $\|f(x)-T(x)\|<\epsilon$ for all $x\in E_1$.}

In 1978, Th. M. Rassias $\cite{RAS1}$ established a generalization of the Hyers' result as the first theorem in the subject of stability of functional equations which allows the Cauchy difference $f(x+y)-f(x)-f(y)$ to be unbounded. This phenomenon has extensively influenced the development of what called Hyers-Ulam-Rassias stability; cf. $\cite{RAS2}, \cite{RAS3}, \cite{RAS4}, \cite{RAS5}$ and $\cite{R-T}$.

During the last decades the problem of Hyers--Ulam--Rassias stability for various functional equations has been widely investigated by many mathematicians. Four methods are used to establish the stability: the Hyers--Ulam sequences, fixed points, invariant means, and sandwich theorems. For a comprehensive account on the stability, the reader is refered to $\cite{CZE1}, \cite{CZE2}, \cite{H-I-R}$. 

In this paper, using Hyers sequence $\cite{HYE}$ and some ideas of $\cite{PAR1}$ and $\cite{RAS1}$ we study the Pexiderized stability of $n$-cocycles and $n$-coboundaries and investigate approximately vanishing of topological cohomology groups as well. In particular, for $n=1$, our results can be regarded as  generalizations of C.-G. Park's results on derivations $\cite{PAR2}$ and multilinear mappings $\cite{PAR1}$.

Throughout this paper, all spaces are assumed to be over the complex field $\C$.

\section {Stability of Cocycles and Coboundaries.}

Throughout this section, $A$ denotes a normed algebra and $X$ is a Banach $A$-bimodule. Suppose that $f_1, f_2, f_3:\displaystyle{\prod_{j=1}^n}A\to X$ are mappings. Fix $n\geq 1$ and scalars $\lambda_1,\cdots,\lambda_n$. Set
\begin{eqnarray*}
&&D^n_{\lambda_1,\cdots,\lambda_n}[f_1,f_2,f_3](a_1,b_1,\cdots,a_n,b_n)\\
&&:=\sum_{j=1}^n(f_1(a_1,\cdots,a_{j-1},\lambda_ja_j
+\lambda_jb_j,a_{j+1},\cdots,a_n)\\
&&-\lambda_jf_2(a_1,\cdots,a_{j-1},a_j,a_{j+1},\cdots,a_n)-\lambda_jf_3(a_1,\cdots,a_{j-1},b_j,a_{j+1},\cdots,a_n)),
\end{eqnarray*}
and
\begin{eqnarray*}
\delta^0x(a):=ax-xa,\\
&&\delta^n[f_1, f_2, f_3](a_1,a_2,\cdots,a_{n+1})\\
&&:=a_1f_1(a_2,\cdots,a_{n+1})
+\sum_ {j=1}^n (-1)^jf_2(a_1,\cdots,a_{j-1},a_ja_{j+1},a_{j+2},\cdots,a_{n+1})\\
&&+(-1)^{n+1}f_3(a_1,\cdots,a_n)a_{n+1}
\end{eqnarray*}
where $x\in X$ and $a_1,\cdots,a_n,a_{n+1},b_1,\cdots,b_n\in A$.

If $f_1=f_2=f_3=f$ we denote
$D^n_{\lambda_1,\cdots,\lambda_n}[f_1,f_2,f_3]$ and
$\delta^n[f_1, f_2, f_3]$ simply by
$D^n_{\lambda_1,\cdots,\lambda_n}f$ and $\delta^nf$
respectively. A mapping $f:\displaystyle{\prod_{j=1}^n}A\to X$ is called multi-linear (multi-additive) if $D^n_{\lambda_1,\cdots,\lambda_n}f=0$ for all $\lambda_1,\cdots,\lambda_n$ ($D^n_{1,\cdots,1}f=0$).
A multi-linear mapping $f$ is said to be $n$-cocycle
if $\delta^nf=0$. By an $n$-coboundary we mean a mapping of the form $\delta^0(x)$ or $\delta^{n-1}g$ in which $g$ is multi-linear.

\begin{t.} Let $\alpha, \beta$ be positive numbers, $n\geq 1, f_1,f_2,f_3:\displaystyle{\prod_{j=1}^n}A\to X$ be mappings such that
\begin{eqnarray}
\|D^n_{\lambda_1,\cdots,\lambda_n}[f_1,f_2,f_3](a_1,b_1,\cdots,a_n,b_n)\|\leq \alpha
\end{eqnarray}
\begin{eqnarray}
\|\delta^n[f_1,f_2,f_3](a_1,a_2,\cdots,a_{n+1})\|\leq \beta
\end{eqnarray}
for all $a_1,\cdots,a_n,a_{n+1},b_1,\cdots,b_n\in A$ and all $\lambda_1,\cdots,\lambda_n\in\C$.\\
Suppose that for each $1\leq k\leq 3$, $f_k(a_1,\cdots,a_n)$ vanishes if $a_i=0$ for any $i$. Then there exists a unique $n$-cocycle $F$ such that
\begin{eqnarray*}
\|f_1(a_1,\cdots,a_n)-F(a_1,\cdots,a_n)\|\leq 3.2^n\alpha\\
\|f_2(a_1,\cdots,a_n)-F(a_1,\cdots,a_n)\|\leq 3(1+\frac{1}{n}).2^n\alpha\\
\|f_3(a_1,\cdots,a_n)-F(a_1,\cdots,a_n)\|\leq 6.2^n\alpha\\
\end{eqnarray*}
Furthermore, if $f_1$ is continuous at a point $(e_1,\cdots,e_n)$ of $\displaystyle{\prod_{j=1}^n}A$ then $F$ is continuous on whole $\displaystyle{\prod_{j=1}^n}A$.\end{t.}

{\bf Proof.} We shall establish the theorem in three steps:

{\bf Step (I): Existence of the multi-linear mapping $F$}

Let $1\leq i\leq n$ be fixed. Putting $\lambda_1=\cdots=\lambda_n=1,
b_1=\cdots=b_n=0$ in $(1)$ we get
\begin{eqnarray*}
\|\sum_{j=1}^nf_1(a_1,\cdots,a_j,\cdots,a_n)-f_2(a_1,\cdots,a_j,\cdots,a_n)\|\leq \alpha,
\end{eqnarray*}
whence
\begin{eqnarray}
\|f_1(a_1,\cdots,a_n)-f_2(a_1,\cdots,a_n)\|\leq \frac{\alpha}{n}
\end{eqnarray}
for all $a_1,\cdots,a_n\in A$. Putting $a_j=(1-\delta_{ij})b_j$ in $(1)$ we obtain
\begin{eqnarray}
\|f_1(b_1,\cdots,b_n)-f_3(b_1,\cdots,b_n)\|\leq \alpha
\end{eqnarray}
for all $b_1,\cdots,b_n\in A$.
Putting $\lambda_1=\cdots=\lambda_n=1, b_j=\delta_{ij}a_i$ in $(1)$ we get
\begin{eqnarray*}
&&\|\sum_{j\in \{1,\cdots,n\}-\{i\}}(f_1(a_1,\cdots,a_j,\cdots,a_n)-f_2(a_1,\cdots,a_j,\cdots,a_n))\\
&&+f_1(a_1,\cdots,a_{i-1},2a_i,a_{i+1},\cdots,a_n)-f_2(a_1,\cdots,a_n)-f_3(a_1,\cdots,a_n)\|\leq \alpha
\end{eqnarray*}
so that
\begin{eqnarray*}
&&\|f_1(a_1,\cdots,a_{i-1},2a_i,a_{i+1},\cdots,a_n)-2f_1(a_1,\cdots,a_n)\|\\
&&\leq \|\sum_{j\in \{1,\cdots,n\}-\{i\}}f_1(a_1,\cdots,a_j,\cdots,a_n)-f_2(a_1,\cdots,a_j,\cdots,a_n)\\
&&+f_1(a_1,\cdots,a_{i-1},2a_i,a_{i+1},\cdots,a_n)-f_2(a_1,\cdots,a_n)-\\
&&f_3(a_1,\cdots,a_n)\|+\|\sum_{j=1}^nf_2(a_1,\cdots,a_j,\cdots,a_n)-f_1(a_1,\cdots,a_j,\cdots,a_n)\|\\
&&+\|f_3(a_1,\cdots,a_n)-f_1(a_1,\cdots,a_n)\|\leq \alpha+\alpha+\alpha=3\alpha.
\end{eqnarray*}
Hence
\begin{eqnarray}
\|f_1(a_1,\cdots,a_{i-1},2a_i,a_{i+1},\cdots,a_n)-2f_1(a_1,\cdots,a_n)\|\leq3\alpha
\end{eqnarray}
Replacing $a_1,\cdots,a_{i-1}$ by $2a_1,\cdots,2a_{i-1}$, respectively, in $(5)$ we get
\begin{eqnarray*}
\|\frac{1}{2^{i-1}}f_1(2a_1,\cdots,2a_{i-1},a_i,a_{i+1},\cdots,a_n)-\frac{1}{2^i}f_1(2a_1,\cdots,2a_{i-1},2a_i,a_{i+1},\cdots,a_n)\|\leq\frac{3}{2^i}\alpha
\end{eqnarray*}
so that
\begin{eqnarray}
&&\|f_1(a_1,\cdots,a_n)-\frac{1}{2^n}f_1(2a_1,\cdots,2a_n)\|\leq\nonumber\\
&&\sum_{i=1}^n \|\frac{1}{2^{i-1}}f_1(2a_1,\cdots,2a_{i-1},a_i,a_{i+1},\cdots,a_n)-\frac{1}{2^i}f_1(2a_1,\cdots,2a_{i-1},2a_i,a_{i+1},\cdots,a_n)\|\nonumber\\
&&\leq\frac{2^n-1}{2}3\alpha
\end{eqnarray}
Replacing $a_1,\cdots,a_n$ by $2^ja_1,\cdots,2^ja_n$ in $(6)$ we get
\begin{eqnarray*}
\|f_1(2^ja_1,\cdots,2^ja_n)-\frac{1}{2^n}f_1(2^{j+1}a_1,\cdots,2^{j+1}a_n)\|\leq\frac{2^n-1}{2}3\alpha
\end{eqnarray*}
whence
\begin{eqnarray*}
&&\|f_1(a_1,\cdots,a_n)-\frac{1}{2^{mn}}f_1(2^ma_1,\cdots,2^ma_n)\\
&&\|\leq\sum_{j=0}^{m-1} \|\frac{1}{2^{nj}}f_1(2^ja_1,\cdots,2^ja_n)-\frac{1}{2^{n+nj}}f_1(2^{j+1}a_1,\cdots,2^{j+1}a_n)\|\\
&&\leq\frac{2^n-1}{2}3\alpha\sum_{j=0}^{m-1}\frac{1}{2^{nj}}.
\end{eqnarray*}
Hence
\begin{eqnarray}
\|f_1(a_1,\cdots,a_n)-\frac{1}{2^{mn}}f_1(2^ma_1,\cdots,2^ma_n)\|\leq
3(1-\frac{1}{2^{mn}})2^n\alpha
\end{eqnarray}
for all $m$ and all $a_1,\cdots,a_n\in A$. Furthermore,
\begin{eqnarray}
\|\frac{1}{2^{m_1n}}f_1(2^{m_1}a_1,\cdots,2^{m_1}a_n)-\frac{1}{2^{m_2n}}f_1(2^{m_2}a_1,\cdots,2^{m_2}a_n)\|\leq\nonumber\\
\frac{2^n-1}{2}3\alpha\sum_{j=m_1}^{m_2-1}(\frac{1}{2^n})^j
\end{eqnarray}
for all $m_2>m_1$.

Inequality $(8)$ shows that the sequence $\{\frac{1}{2^{mn}}f_1(2^ma_1,\cdots,2^ma_n)\}$ is Cauchy in the Banach module $X$ and so is convergent. Set
\begin{eqnarray}
F(a_1,\cdots,a_n):=\displaystyle{\lim_{m\to\infty}}\frac{1}{2^{mn}}f_1(2^ma_1,\cdots,2^ma_n).
\end{eqnarray}
Inequality $(7)$ yields
\begin{eqnarray*}
\|f_1(a_1,\cdots,a_n)-F(a_1,\cdots,a_n)\|\leq 3.2^n\alpha
\end{eqnarray*}
By $(3)$,
\begin{eqnarray*}
\|2^{-mn}f_1(2^ma_1,\cdots,2^ma_n)-2^{-mn}f_2(2^ma_1,\cdots,2^ma_n)\|\leq \frac{\alpha}{2^{mn}n}.
\end{eqnarray*}
Using $(9)$ we have
\begin{eqnarray}
F(a_1,\cdots,a_n)=\displaystyle{\lim_{m\to\infty}}\frac{1}{2^{mn}}f_2(2^ma_1,\cdots,2^ma_n).
\end{eqnarray}
By $(5)$ and $(7)$ we get
\begin{eqnarray*}
&&\|f_2(a_1,\cdots,a_{i-1},2a_i,a_{i+1},\cdots,a_n)-2f_2(a_1,\cdots,a_n)\|\\
&&\leq\|f_2(a_1,\cdots,a_{i-1},2a_i,a_{i+1},\cdots,a_n)-f_1(a_1,\cdots,a_{i-1},2a_i,a_{i+1},\cdots,a_n)\\
&&+\|f_1(a_1,\cdots,a_{i-1},2a_i,a_{i+1},\cdots,a_n)-2f_1(a_1,\cdots,a_n)\|\\
&&+2\|f_1(a_1,\cdots,a_n)-f_2(a_1,\cdots,a_n)\|\\
&&\leq\frac{\alpha}{n}+3\alpha+2\frac{\alpha}{n}\\
&&=3(1+\frac{1}{n})\alpha
\end{eqnarray*}
so that
\begin{eqnarray}
\|f_2(a_1,\cdots,a_{i-1},2a_i,a_{i+1},\cdots,a_n)-2f_2(a_1,\cdots,a_n)\|\leq 3(1+\frac{1}{n})\alpha.
\end{eqnarray}
As the same way as we obtained inequality $(7)$, one can deduce from $(11)$ that
\begin{eqnarray*}
\|f_2(a_1,\cdots,a_n)-\frac{1}{2^{mn}}f_2(2^ma_1,\cdots,2^ma_n)\|\leq
3(1+\frac{1}{n})(1-\frac{1}{2^{mn}})2^n\alpha.
\end{eqnarray*}
Letting $m$ tend to $\infty$ we obtain
\begin{eqnarray*}
\|f_2(a_1,\cdots,a_n)-F(a_1,\cdots,a_n)\|\leq 3(1+\frac{1}{n}).2^n\alpha
\end{eqnarray*}
Similarly, by applying $(4)$ we obtain
\begin{eqnarray}
F(a_1,\cdots,a_n)=\displaystyle{\lim_{m\to\infty}}\frac{1}{2^{mn}}f_3(2^ma_1,\cdots,2^ma_n),
\end{eqnarray}
and
\begin{eqnarray*}
\|f_3(a_1,\cdots,a_n)-F(a_1,\cdots,a_n)\|\leq 6.2^n\alpha.
\end{eqnarray*}

Replacing $a_1,\cdots,a_n,b_1,\cdots,b_n,\lambda_1,\cdots,\lambda_i,\cdots,\lambda_n$ by\\
$2^ma_1,\cdots,2^ma_n,0,\cdots,0,1,\cdots,\lambda,\cdots,1$, respectively, in $(1)$ we get
\begin{eqnarray*}
&&\frac{1}{2^{mn}}D^n_{1,\cdots,\lambda_i,\cdots,1}[f_1,f_2,f_3](2^ma_1,0,\cdots,2^ma_{i-1},0,2^ma_i,2^mb_i,2^ma_{i+1},0,\cdots,2^ma_n,0)\|\\
&&\leq \frac{\alpha}{2^{mn}}
\end{eqnarray*}
Passing to the limit as $m\to\infty$ we conclude that
\begin{eqnarray*}
&&F(a_1,\cdots,a_{i-1},\lambda_ia_i+\lambda_ib_i,a_{i+1},\cdots,a_n)\\
&&=\lambda_iF(a_1,\cdots,a_{i-1},a_i,a_{i+1},\cdots,a_n)+\lambda_iF(a_1,\cdots,a_{i-1},b_i,a_{i+1},\cdots,a_n).
\end{eqnarray*}
Therefore $F$ is linear in the $i$-th variable for each $i=1,\cdots,n$.

If $F':\displaystyle{\prod_{j=1}^n}A\to X$ is a multi-linear mapping with $\|f(a_1,\cdots,a_n)-F'(a_1,\cdots,a_n)\|\leq 3.2^n\alpha$ for all $a_1,\cdots,a_n\in A$ then
\begin{eqnarray*}
&&\|F(a_1,\cdots,a_n)-F'(a_1,\cdots,a_n)\|\\
&&=\lim_{m\to\infty}2^{-mn}\|f(2^ma_1,\cdots,2^ma_n)-F'(2^ma_1,\cdots,2^ma_n)\|\\
&&\leq\lim_{m\to\infty}\frac{3.2^n\alpha}{2^{mn}}\\
&&=0
\end{eqnarray*}
whence $F=F'$.

{\bf Step (II): Proving $F$ to be cocycle.}

For each fixed $1\leq i\leq n$ and $a_1,\cdots,a_n\in A$ one can apply $(11)$ and induction on $m$ to prove
\begin{eqnarray}
\|2^{-m}f_2(a_1,\cdots,a_{i-1},2^ma_i,a_{i+1},\cdots,a_n)-f_2(a_1,\cdots,a_{i-1},a_i,a_{i+1},\cdots,a_n)\|\nonumber\\
\leq 3(1-2^{-m})(1+\frac{1}{n})\alpha
\end{eqnarray}
Now we can replace $a_i$ by $2^ma_i$ in $(13)$ to get
\begin{eqnarray}
&&\|2^{-(n+1)m}f_2(2^ma_1,\cdots,2^ma_{i-1},2^{2m}a_i,2^ma_{i+1},\cdots,2^ma_n)\\
&&-2^{-mn}f_2(2^ma_1,\cdots,2^ma_{i-1},2^ma_i,2^ma_{i+1},\cdots,2^ma_n)\|\leq 3(\frac{1}{2^{mn}}-\frac{1}{2^{m(n+1)}})\alpha
\end{eqnarray}
Then $(10)$ and $(14)$ yield
\begin{eqnarray}
F(a_1,\cdots,a_n)=\displaystyle{\lim_{m\to\infty}}\frac{1}{2^{m(n+1)}}f_2(2^ma_1,\cdots,2^ma_{i-1},2^{2m}a_i,2^ma_{i+1},\cdots,2^ma_n)
\end{eqnarray}
By $(2)$, we have
\begin{eqnarray*}
&&\|2^{-(n+1)m}\delta^n[f_1,f_2,f_3](a_1,\cdots,a_{n+1})\|=\|2^{-mn}a_1f_1(2^ma_2,\cdots,2^ma_{n+1})\\
&&+2^{-(n+1)m}\sum_ {j=1}^n (-1)^jf_2(2^ma_1,\cdots,2^ma_{j-1},2^{2m}a_ja_{j+1},a_{j+2},\cdots,2^ma_{n+1})\\
&&+(-1)^{n+1}2^{-mn}f_3(2^ma_1,\cdots,2^ma_n)a_{n+1}\|\\
&&\leq 2^{-(n+1)m}\beta
\end{eqnarray*}
for all $m$ and all $a_1,\cdots,a_{n+1}\in A$.\\

Next by passing to the limit as $m\to\infty$ and noting $(9)$, $(12)$ and $(15)$ we get
\begin{eqnarray*}
\delta^nF(a_1,\cdots,a_{n+1})&=&a_1F(a_2,\cdots,a_{n+1})\\
&+&\sum_ {j=1}^n (-1)^jF(a_1,\cdots,a_{j-1},a_ja_{j+1},a_{j+2},\cdots,a_{n+1})\\
&+&(-1)^{n+1}F(a_1,\cdots,a_n)a_{n+1}=0.
\end{eqnarray*}
for all $a_1,\cdots,a_{n+1}\in A$. Hence $F$ is a cocycle.

{\bf Step (III): Continuity of $F$}

We use the strategy of Hyers $\cite{HYE}$. If $F$ were not continuous at the point $(e_1,\cdots,e_n)$ then there would be an integer $P$ and a sequence $\{(a_1^m,\cdots,a_n^m)\}$ of $\displaystyle{\prod_{j=1}^n}A$ converging to zero such that $\|F(a_1^m,\cdots,a_n^m)\|>\frac{1}{P}$. Let $K$ be an integer greater than $7P2^n\alpha$. Since $\displaystyle{\lim_{m\to\infty}}f_1(K(a_1^m,\cdots,a_n^m)+(e_1,\cdots,e_n))=f_1(e_1,\cdots,e_n)$, there is an integer $N$ such that $\|f_1(K(a_1^m,\cdots,a_n^m)+(e_1,\cdots,e_n))-f_1(e_1,\cdots,e_n)\|<2^n\alpha$ for all $n\geq N$. Hence 
\begin{eqnarray*}
7.2^n\alpha&<&\frac{K}{P}<\|F(K(a_1^m,\cdots,a_n^m))\|\\
&=&\|F(K(a_1^m,\cdots,a_n^m)+(e_1,\cdots,e_n))-F(e_1,\cdots,e_n)\|\\
&\leq&\|F(K(a_1^m,\cdots,a_n^m)+(e_1,\cdots,e_n))-f_1(K(a_1^m,\cdots,a_n^m)-(e_1,\cdots,e_n))\|\\
&&+\|f_1(K(a_1^m,\cdots,a_n^m)-(e_1,\cdots,e_n))-f_1(e_1,\cdots,e_n)\|\\
&&+\|f_1(e_1,\cdots,e_n)-F(e_1,\cdots,e_n)\|\\
&&\leq 3.2^n\alpha+2^n\alpha+3.2^n\alpha\\
&&=7.2^n\alpha
\end{eqnarray*}
for all $n>N$, a contradiction. Now the multi-linearity of $F$ guarantee continuity of $F$ on  whole $\displaystyle{\prod_{j=1}^n}A.\Box$

\begin{t.} Let $\alpha, \beta, \gamma$ be positive numbers, $x\in X$ and $f_1,f_2,f_3:A\to X$ be  mappings such that
\begin{eqnarray*}
\|D^1_{\lambda}[f_1,f_2,f_3](a,b)\|\leq \alpha
\end{eqnarray*}
\begin{eqnarray*}
\|\delta^1[f_1,f_2,f_3](a,b)\|\leq \beta
\end{eqnarray*}
\begin{eqnarray}
\|ax-xa-f_1(a)\|\leq\gamma
\end{eqnarray}
for all $a,b\in A$ and all $\lambda\in\C$.\\
Suppose that for each $1\leq k\leq 3$, $f_k(0)=0$. Then there exists a $1$-cocycle $F$ such that
\begin{eqnarray*}
\|f_1(a)-F(a)\|\leq 6\alpha\\
\|f_2(a)-F(a)\|\leq 12\alpha\\
\|f_3(a)-F(a)\|\leq 12\alpha\\
F(a)=ax-xa
\end{eqnarray*}
for all $a\in A$.\end{t.}

{\bf Proof.} By Theorem 2.1, there is a unique $1$-cocycle $F$
defined by $F(a):=\displaystyle{\lim_{m\to\infty}}2^{-m}f_1(2^ma)$
satisfying the required inequalities. It follows from $(16)$ that
$\|ax-xa-2^{-m}f_1(2^ma)\|\leq 2^{-m}\gamma$. Passing to the limit we
conclude that $F(a)=ax-xa.\Box$

\begin{r.} {\rm  Theorem 2.2. gives us the Hyers--Ulam stability of any one of the following function equations:\\
(i) $f(ab)=af(b)+f(a)b$; cf. $\cite{PAR2}$\\
(ii) $af(b)=f(a)b$\\
(iii) $f(ab)=af(b)$\\
(iv) $f(ab)=f(a)b$\\
together with the Cauchy equation $f(a+b)=f(a)+f(b)$.\\
To see this use Theorem 2.2. with $f_1=f_2=f_3=f$ to get (i); $f_1=f_3=f, f_2=0$ to obtain (ii); $f_1=f_2=f, f_3=0$ to get (iii); and $f_1=0, f_2=f_3=f$ to obtain (iv).}\end{r.}

\begin{p.} Let $A$ be linearly spanned by a set $S\subseteq A, \alpha, \beta$ be positive numbers, $n\geq 1, f_1,f_2,f_3:\displaystyle{\prod_{j=1}^n}A\to X$ be mappings such that
\begin{eqnarray*}
\|D^n_{\lambda_1,\cdots,\lambda_n}[f_1,f_2,f_3](a_1,b_1,\cdots,a_n,b_n)\|\leq \alpha
\end{eqnarray*}
for all $a_1,\cdots,a_n,b_1,\cdots,b_n\in A$ and all $\lambda_1,\cdots,\lambda_n\in \T=\{z\in\C: |z|=1\}$; and
\begin{eqnarray*}
\|\delta^n[f_1,f_2,f_3](a_1,a_2,\cdots,a_{n+1})\|\leq \beta
\end{eqnarray*}
for all $a_1,\cdots,a_{n+1}\in S$.\\
Suppose that for each $1\leq k\leq 3$, $f_k(a_1,\cdots,a_n)$ vanishes if $a_i=0$ for any $i$. Then there exists a unique $n$-cocycle $F$ such that
\begin{eqnarray*}
\|f_1(a_1,\cdots,a_n)-F(a_1,\cdots,a_n)\|\leq 3.2^n\alpha\\
\|f_2(a_1,\cdots,a_n)-F(a_1,\cdots,a_n)\|\leq 3(1+\frac{1}{n}).2^n\alpha\\
\|f_3(a_1,\cdots,a_n)-F(a_1,\cdots,a_n)\|\leq 6.2^n\alpha\\
\end{eqnarray*}\end{p.}

{\bf Proof.} By the same argument as in the proof of Theorem 2.1 there exists a unique multi-additive mapping $F$ satisfying the required inequalities such that $\delta^nF(a_1,a_2,\cdots,a_{n+1})=0$ holds for all $a_1,\cdots,a_{n+1}\in S$. 

Fix $1\leq i\leq n$. Assume that $\lambda_i\in \C$ and $\lambda\neq 0$. If $N$ is a positive integer greater than $4\|a_i\|$ then $\|\frac{a_i}{N}\|<\frac{1}{4}<1-\frac{2}{3}=\frac{1}{3}$. By Theorem 1 of $\cite{K-P}$ there are three numbers $z_1,z_2,z_3\in \T$ such that $3\frac{a_i}{N}=z_1+z_2+z_3$. By virtue of the multi-additivity of $F$ we easily conclude that $F$ is multi-linear. Since each element of $A$ is a linear combination of elements of $S$, we infer that $F$ is a cocycle.$\Box$

\begin{p.}  Let $A$ be linearly spanned by a set $S\subseteq A, \alpha, \beta$ be positive numbers, $n\geq 1, f_1,f_2,f_3:\displaystyle{\prod_{j=1}^n}A\to X$ be mappings such that
\begin{eqnarray*}
\|D^n_{\lambda_1,\cdots,\lambda_n}[f_1,f_2,f_3](a_1,b_1,\cdots,a_n,b_n)\|\leq \alpha
\end{eqnarray*}
for all $a_1,\cdots,a_n,b_1,\cdots,b_n\in A$ and all $\lambda_1,\cdots,\lambda_n\in\{1, {\rm i}\}$ and
\begin{eqnarray*}
\|\delta^n[f_1,f_2,f_3](a_1,a_2,\cdots,a_{n+1})\|\leq \beta
\end{eqnarray*}
for all $a_1,\cdots,a_{n+1}\in S$.\\
Suppose that for each $1\leq k\leq 3$, $f_k(a_1,\cdots,a_n)$ vanish if $a_i=0$ for any $i$. Assume that for each $1\leq i\leq n$ and each fixed $(a_1,\cdots,a_n)$ the function $t\mapsto f(a_1,\cdots,a_{i-1},ta_i,a_{i+1},\cdots,a_n)$ is continuous on $\R$. Then there exists a unique $n$-cocycle $F$ such that
\begin{eqnarray*}
\|f_1(a_1,\cdots,a_n)-F(a_1,\cdots,a_n)\|\leq 3.2^n\alpha\\
\|f_2(a_1,\cdots,a_n)-F(a_1,\cdots,a_n)\|\leq 3(1+\frac{1}{n}).2^n\alpha\\
\|f_3(a_1,\cdots,a_n)-F(a_1,\cdots,a_n)\|\leq 6.2^n\alpha\\
\end{eqnarray*}\end{p.}

{\bf Proof.} By the same argument as in the proof of Theorem 2.1 there exists a unique multi-additive mapping $F$ satisfying the required inequalities such that $\delta^nF(a_1,a_2,\cdots,a_{n+1})=0$ holds for all $a_1,\cdots,a_{n+1}\in S$. 

Fix $1\leq i\leq n$. By the same reasoning as in the proof of theorem of $\cite{RAS1}$, the mapping $F$ is multi-$\R$-linear. Since $\C$ as a vector space over $\R$ is generated by $\{1,{\rm i}\}$, we conclude that $F$ is multi-$\C$-linear. Since each element of $A$ is a linear combination of elements of $S$, we infer that $F$ is a cocycle.$\Box$

\begin{t.}  Let $\alpha, \beta, \gamma, \eta$ be positive numbers, $n\geq 2, f_1,f_2,f_3:\displaystyle{\prod_{j=1}^n}A\to X$ and $g_1,g_2,g_3:\displaystyle{\prod_{j=1}^{n-1}}A\to X$ be  mappings such that
\begin{eqnarray*}
&&\|D^n_{\lambda_1,\cdots,\lambda_n}[f_1,f_2,f_3](a_1,b_1,\cdots,a_n,b_n)\|\leq \alpha,\\
&&\|\delta^n[f_1,f_2,f_3](a_1,a_2,\cdots,a_{n+1})\|\leq \beta,\\
&&\|D^{n-1}_{\lambda_1,\cdots,\lambda_{n-1}}[g_1,g_2,g_3](a_1,b_1,\cdots,a_{n-1},b_{n-1})\|\leq \gamma
\end{eqnarray*}
\begin{eqnarray}
&&\|\delta^{n-1}[g_1,g_2,g_3](a_1,a_2,\cdots,a_n)-f_1(a_1,\cdots,a_{n+1})\|\leq\eta
\end{eqnarray}
for all $a_1,\cdots,a_{n-1},a_n,a_{n+1},b_1,\cdots,b_{n-1},b_n\in A$ and all $\lambda_1,\cdots,\lambda_{n-1},\lambda_n\in\C$.\\
Suppose that for each $1\leq k\leq 3$, $f_k(a_1,\cdots,a_n)$ and
$g_k(a_1,\cdots,a_{n-1})$ vanish if $a_i=0$ for any $i$ and $g_1$ is continuous at a point of $\displaystyle{\prod_{j=1}^{n-1}}A$. Then
there exist a unique $n$-cocycle $F$ and a unique continuous multi-linear
mapping $G:\displaystyle{\prod_{j=1}^{n-1}}A\to X$ such that
\begin{eqnarray*}
&&\|f_1(a_1,\cdots,a_n)-F(a_1,\cdots,a_n)\|\leq 3.2^n\alpha,\\
&&\|f_2(a_1,\cdots,a_n)-F(a_1,\cdots,a_n)\|\leq 3(1+\frac{1}{n}).2^n\alpha,\\
&&\|f_3(a_1,\cdots,a_n)-F(a_1,\cdots,a_n)\|\leq 6.2^n\alpha,\\
&&\|g_1(a_1,\cdots,a_{n-1})-G(a_1,\cdots,a_{n-1})\|\leq 3.2^n\gamma,\\
&&\|g_2(a_1,\cdots,a_{n-1})-G(a_1,\cdots,a_{n-1})\|\leq 3(1+\frac{1}{n}).2^n\gamma,\\
&&\|g_3(a_1,\cdots,a_{n-1})-G(a_1,\cdots,a_{n-1})\|\leq 6.2^n\gamma,
\end{eqnarray*}
and
\begin{eqnarray*}
F=\delta^nG.
\end{eqnarray*}\end{t.}

{\bf Proof.} Theorem 2.1 gives rise the existence of a unique $n$-cocycle $F$
with requested properties. Using the same reasoning as in the
proof of Theorem 2.1 one can show that there exists a unique continuous
multi-linear mapping $G$ defined by
\begin{eqnarray}
G(a_1,\cdots,a_{n-1}):=\displaystyle{\lim_{m\to\infty}}\frac{1}{2^{m(n-1)}}g_1(2^ma_1,\cdots,2^ma_{n-1})
\end{eqnarray}
satisfying the required inequalities and
\begin{eqnarray}
&&G(a_1,\cdots,a_{n-1})=\displaystyle{\lim_{m\to\infty}}\frac{1}{2^{m(n-1)}}g_3(2^ma_1,\cdots,2^ma_{n-1})\\
&&G(a_1,\cdots,a_{n-1})=\displaystyle{\lim_{m\to\infty}}\frac{1}{2^{mn}}g_2(2^ma_1,\cdots,2^ma_{i-1},2^{2m}a_i,2^ma_{i+1},\cdots,2^ma_{n-1})
\end{eqnarray}
Clearly,
\begin{eqnarray*}
\delta^{n-1}G(a_1,\cdots,a_n)=\displaystyle{\lim_{m\to\infty}}2^{-mn}\delta^{n-1}[g_1,g_2,g_3](2^ma_1,\cdots,2^ma_n).
\end{eqnarray*}
Inequality $(17)$ yields
\begin{eqnarray*}
&&\|2^{-mn}\delta^{n-1}[g_1,g_2,g_3](2^ma_1,\cdots,2^ma_n)-2^{-mn}f_1(2^ma_1,\cdots,2^ma_n)\|=\\
&&\|2^{-m(n-1)}a_1g_1(2^ma_2,\cdots,2^ma_{n+1})\\
&&+2^{-nm}\sum_{j=1}^{n-1}(-1)^jg_2(2^ma_1,\cdots,2^ma_{j-1},2^{2m}a_ja_{j+1},a_{j+2},\cdots,2^ma_n)\\
&&+(-1)^n2^{-m(n-1)}g_3(2^ma_1,\cdots,2^ma_{n-1})a_n-2^{-mn}f_1(2^ma_1,\cdots,2^ma_n)\|\leq 2^{-mn}\eta
\end{eqnarray*}
Letting $m\to\infty$ and using $(9)$, $(18)$, $(19)$ and $(20)$ we conclude that
\begin{eqnarray*}
\|a_1G(a_2,\cdots,a_n)\\
&+&\sum_ {j=1}^{n-1} (-1)^jG(a_1,\cdots,a_ja_{j+1},a_{j+2},\cdots,a_n)\\
&+&(-1)^nG(a_1,\cdots,a_{n-1})a_n-F(a_1,\cdots,a_n)\|=0
\end{eqnarray*}
Thus $\delta^{n-1}(G)=F.\Box$

\begin{r.} {\rm There are statements similar to Propositions 2.4, 2.5 for coboundaries.}\end{r.}

\section{Vanishing of Cohomology Groups.}

Throughout this section, $A$ denotes a Banach algebra and $X$ is a Banach $A$-bimodule. For $n=0,1,2,\cdots$, let $ C^n(A,X)$ be the Banach space of all
bounded $n$-linear mappings from $A\times\cdots\times A$ into $X$
equipped with multi-linear operator norm $\| f\|=\sup \{\|
f(a_1,\cdots,a_n)\|: a_i\in A, \|a_i\|\leq 1, 1\leq i\leq n\}$,
and $C^0(A,X)=X$. The elements of $C^n(A,X)$ are called
$n$-dimensional cochains. Consider the sequence
$$0 \to C^0(A,X)\stackrel{\delta^0}{\to} C^1(A,X)\stackrel{\delta^1}{\to}\cdots~~(\widetilde{C}(A,X))$$

It is not hard to show that the above sequence is a complex, i.e. for each $n, \delta^{n+1}\circ \delta^n=0$; cf. $\cite{RIN}$.

$\widetilde{C}(A,X)$ is called the standard cohomology complex or
Hochschild--Kamowitz complex for $A$ and $X$. The $n$-th
cohomology group of $\widetilde{C}(A,X)$ is said to be
$n$-dimensional (ordinary or Hochschild) cohomology group of $A$
with coefficients in $X$ and denoted by $H^n(A,X)$. The spaces
${\rm Ker}\delta ^n$ and ${\rm Im}\delta^{n-1}$ are denoted by
$Z^n(A,X)$ and $B^n(A,X)$, respectively. Hence
$H^n(A,X)=Z^n(A,X)/B^n(A,X)$. The cohomology groups of small
dimensions $n=0,1,2,3$ are very important and applicable.

$H^0(A,X)=Z^0(A,X)$ is the so-called center of $X$.

Any element of $$Z^1(A,X)=\{d: A \to X; d {\rm ~ is~
bounded~ and~ linear,~ and~} d(ab)=ad(b)+d(a)b\}$$
is called a derivation of $A$ in $X$ and any element of
$B^1(A,X)=\{d_x: A \to X; d_x(a)=ax-xa, a\in A, x\in
X\}$ is called an inner derivation. The Banach algebra $A$ is said to
be contractible if $H^1(A,X)=0$ for all $X$ and to be amenable (according to Johnson) if $H^1(A,X^*)=0$ for all $X$; cf. $\cite{JOH}$.

$H^2(A,X)$ is the equivalence classes of singular extensions of $A$ by $X$; cf. $\cite{B-D-L}$.

$H^3(A,X)$ can be used in the study of stable properties of Banach algebras; cf. $\cite{JAR}$.

For $n\geq 4$ there is no known interesting interpretation of $H^n(A,X)$. But their vanishing is what homological dimension is about $\cite{HEL2}$.
Given $n\geq 1$, by an {\it approximate $n$-cocycle} we mean a mapping $f:\displaystyle{\prod_{j=1}^n}A\to X$ which is continuous at a point and $f(a_1,\cdots,a_n)=0$ whenever $a_i=0$ for any $i$, and such that 
\begin{eqnarray*}
\|D^n_{\lambda_1,\cdots,\lambda_n}f(a_1,b_1,\cdots,a_n,b_n)\|\leq \alpha
\end{eqnarray*}
\begin{eqnarray*}
\|\delta^nf(a_1,a_2,\cdots,a_{n+1})\|\leq \beta
\end{eqnarray*}
for some positive numbers $\alpha$ and $\beta$ and for all $a_1,\cdots,a_n,a_{n+1},b_1,\cdots,b_n\in A$ and all $\lambda_1,\cdots,\lambda_n\in\C$.

Given $n\geq 2$, by an {\it approximate $n$-coboundary} we mean a mapping of the form $\delta^{n-1}g$ in which $g:\displaystyle{\prod_{j=1}^{n-1}}A\to X$ that is continuous at a point and $g(a_1,\cdots,a_{n-1})=0$ whenever $a_i=0$ for any $i$, and such that 
\begin{eqnarray*}
\|D^{n-1}_{\lambda_1,\cdots,\lambda_{n-1}} g(a_1,b_1,\cdots,a_{n-1},b_{n-1})\|\leq \gamma
\end{eqnarray*}
for some positive number $\gamma$ and for all $a_1,\cdots,a_{n-1},b_1,\cdots,b_{n-1}\in A$ and all $\lambda_1,\cdots,\lambda_{n-1}\in\C$.\\
By an {\it approximate $1$-coboundary} we mean a mapping of the form
$\delta^0(a)=ax-xa$ for some $x\in X$, i.e. a usual $1$-coboundary.

If every approximate $n$-cocycle $f$ is near an approximate
$n$-coboundary, in the sense that there exist $\eta>0$ and an approximate $n$-coboundary $h$ such that $\|h(a_1,\cdots,a_n)-f(a_1,\cdots,a_n)\|\leq \eta$ for all $a_1,\cdots,a_n\in A$, we say the $n$-th cohomology group of $A$ with coefficients in $X$ approximately vanishes.

\begin{t.} For a positive integer $n, H^n(A,X)=0$ if and only if the $n$-th cohomology group of $A$ in $X$ approximately vanishes.\end{t.}

{\bf Proof.} Suppose that $H^n(A,X)=0$ and $f$ is an approximate $n$-cocycle. By Theorem 2.1 there is an $n$-cocycle $F\in Z^n(A,X)$ such that $\|F(a_1,\cdots,a_n)-f(a_1,\cdots,a_n)\|\leq 3.2^n\alpha$ where $\alpha$ is given by $(1)$. Since $H^n(A,X)=0$, there exists $G\in C^{n-1}(A,X)$ such that $\delta^{n-1}G=F$. Hence $\|\delta^{n-1}G(a_1,\cdots,a_n)-f(a_1,\cdots,a_n)\|\leq 3.2^n\alpha$. Hence $f$ is approximated by an {\it approximate coboundary}. 

For the converse, let $F\in Z^n(A,X)$. Then $F$ is trivially an approximate $n$-cocycle. Since $n$-th cohomology group of $A$ in $X$ approximately vanishes, there exist $\eta>0$ and an approximate $n$-coboundary $h$ such that $\|h(a_1,\cdots,a_n)-F(a_1,\cdots,a_n)\|\leq\eta$. By Theorem 2.2 and Theorem 2.6 there exist $G\in C^{n-1}(A,X)$ such that $F=\delta^{n-1}G$. Hence $F\in B^n(A,X)$. Therefore $H^n(A,X)=0.\Box$

\begin{c.} The Banach algebra $A$ is contractible if and only if it is {\it approximately contractible}, i.e every continuous approximate derivation from $A$ into any Banach $A$-bimodule is near an inner derivation.\end{c.} (See $\cite{MOS}$ for another approach)

\begin{c.} The Banach algebra $A$ is amenable if and only if it is {\it approximately amenable}, i.e every continuous approximate derivation from $A$ into a dual Banach $A$-bimodule is near an inner derivation.\end{c.}

{\bf Address.}\\
Mohammad Sal Moslehian\\
Dept. of Math., Ferdowsi Univ., P. O. Box 1159, Mashhad 91775, Iran\\
E-mail: moslehian@ferdowsi.um.ac.ir\\
Home: http://www.um.ac.ir/$\sim$moslehian/
\end{document}